\documentclass[12pt,leqno]{article}
\usepackage[pdftex]{graphicx}
\usepackage{amsfonts,amssymb,amsmath}
\usepackage{ccaption,mathrsfs,amsmath,amssymb,amsthm}
\usepackage{amsmath, amscd, amsfonts,  tikz, times}

\newcounter{conjecture}
\newcounter{remark}

\newcommand{\eqnsection}{
    \renewcommand{\theequation}{\arabic{equation}}
    \makeatletter
    \csname @addtoreset\endcsname{equation}{section}
    \makeatother}
\newtheorem{theorem}{Theorem}

\newtheorem{lemma}{Lemma}

\newcommand{\dd}{\delta}

\newcommand{\lar}{\longrightarrow}

\newcommand{\aaa}{\alpha}

\newcommand{\CC}{\mathbb{C}}

\def \be{\begin{equation}}
\def \ee{\end{equation}}
\def \bt{\begin{theorem}}
\def \et{\end{theorem}}
\def \bea{\begin{eqnarray}}
\def \eea{\end{eqnarray}}
\def \bas{\begin{eqnarray*}}
\def \eas{\end{eqnarray*}}

\newcommand {\rrr}[1]{(\ref{#1})}

\def \ga{\gamma}

\def \Om{\Omega}

\def \si{\sigma}

\def \th{\theta}

\def \ff{\infty}

\def \DD{{\mathbb D}}

\def \HH{{\mathbb H}}

\def \RR{{\mathbb R}}

\def \({\left(}
\def \){\right)}

\def \vski{\vspace{12pt}}

\def \bc{\begin{center} }
\def \ec{\end{center} }
\def \bs{\begin{slide} }
\def \es{\end{slide} }

\def\square{{\vcenter{\vbox{\hrule height.3pt
         \hbox{\vrule width.3pt height5pt \kern5pt
            \vrule width.3pt}
         \hrule height.3pt}}}}

\newcounter{cccases}
\eqnsection
\begin{document}

\title{Many proofs that $\sum_{n=1}^{\ff} \frac{1}{n^2} = \frac{\pi^2}{6}$ can be found in the conformal invariance of planar Brownian motion}

\author{
\begin{tabular}{c}
\textit{Greg Markowsky} \\
Monash University \\
Victoria 3800, Australia \\
gmarkowsky@gmail.com
\end{tabular}}

\bibliographystyle{amsplain}

\maketitle \eqnsection \setlength{\unitlength}{2mm}

\begin{abstract}
\noindent A number of recent new proofs of Euler's celebrated identity are presented, all of which are probabilistic in nature and depend on the conformal invariance of planar Brownian motion.




\end{abstract}

\section*{Introduction}

Leonhard Euler first achieved fame in 1734 with his evaluation of $\sum_{n=1}^{\ff} \frac{1}{n^2}$, a quantity which had fascinated leading mathematicians for roughly a century before him. He arrived at the correct value of $\frac{\pi^2}{6}$ by assuming that the function $\frac{\sin x}{x}$ admitted a factorization into the product of linear terms as though it were a polynomial:

\begin{equation} \label{sineprod}
\frac{\sin x}{x} = (1-\frac{x}{\pi})(1+\frac{x}{\pi})(1-\frac{x}{2\pi})(1+\frac{x}{2\pi})(1-\frac{x}{3\pi})(1+\frac{x}{3\pi}) \ldots
\end{equation}

Armed with this expansion, the result now follows by equating the $x^2$ term in the power series expansion of both sides. However, proving \rrr{sineprod} is not so simple, and in the time of Euler the identity was not rigorously established; it was only proved to complete satisfaction a century later by Weierstrass as part of the theory of infinite products. Since Euler's initial achievement a large number of new proofs have been uncovered; \cite{chappy} contains 14 of them, and there are others scattered throughout the literature.

\vski

Recently, a number of new proofs have been discovered and appear (either explicitly or implicitly) in the papers \cite{meecp}, \cite{megreen}, and \cite{medist}. The purpose of this short paper is to present them explicitly in one place, together with necessary background information. All of these proofs utilize the invariance of Brownian motion under mapping by complex analytic functions; it is hoped that this presentation may serve as a relatively nontechnical introduction for the uninitiated to the charming interplay between complex analysis and probability theory.

\vski

It will be assumed that the reader has a basic familiarity with standard analytic functions and martingale theory. The necessary background in complex analysis can be found in \cite{rud} or \cite{schaum}, and for martingale theory one may consult \cite{fima} or \cite{revyor}; several sources, such as \cite{durBM} and \cite{bassy}, focus on the interaction between the two fields. The main tool in all that follows will be the following major theorem of L\'evy (\cite[Ch. V]{bassy}).

\begin{theorem} \label{holinv}
Let $f$ be analytic and nonconstant on a domain $\Om$, and let $a \in \Om$. Let $B_t$ be a 2-dimensional Brownian motion which starts at $a$ and is stopped at a stopping time $\tau$ such that the set of Brownian paths $\{B_t: 0 \leq t \leq \tau\}$ lie within the closure of $\Om$ a.s. Set

\begin{equation*} \label{}
\si(t) = \int_{0}^{t \wedge \tau} |f'(B_s)|^2 ds.
\end{equation*}

$\si(t)$ is a.s. strictly increasing and continuous, and we let $C_{t} = \si^{-1}(t)$. Then $\hat B_t = f(B_{C_t})$ is a Brownian motion stopped at the stopping time $\si(\tau)$.
\end{theorem}

We will see that many results, including the promised proofs of Euler's sum, flow easily from this theorem paired with a few well-chosen analytic functions. In what follows, $B_t$ will always denote planar Brownian motion, and we will use the notation of the complex plane; for example, we will write $B_t = Re(B_t) + i Im(B_t)$ when necessary, where $Re(B_t)$ and $Im(B_t)$ are 1-dimensional Brownian motions. We will also use the standard notation $E_a$ to denote expectation conditioned on $B_0 = a$ almost surely.

\vski

As a final disclaimer before beginning the proofs, it should be mentioned that very few, if any, of the ideas presented below are truly original to \cite{meecp}, \cite{megreen}, \cite{medist}; many of the same methods have surfaced over the years in various contexts, albeit often packaged somewhat differently. For example, the idea of wrapping the real axis around the circle as a mechanism for deducing summation identities (Proof 2 below) is essentially the same as is used in the Poisson summation formula (see for example \cite{stein}), and has been used also to deduce Jacobi's theta function identity from properties of one-dimensional Brownian motion wound around the circle (\cite{biane}). The distribution of the exit time of one-dimensional Brownian motion from an interval (which is the same as the exit time of of planar Brownian motion from a strip) has been used to deduce Leibniz's formula and related sums in \cite{pityor}, and in Proof 3 below we deduce the same formula from the exit distribution of the strip. In \cite{pace} and in \cite{bourg} (see also \cite{yano}), Euler's formula is deduced from properties of the Cauchy distribution, which has an intimate connection with planar Brownian motion as it is the exit distribution of a half-plane; in fact, this connection is utilized in \cite{bourg}. Expressing Green's function as an infinite product in certain domains had previously been discussed, purely as an analytic method, in \cite{melgreen} and \cite{melgreen2}. And so forth.

\section*{Proof 1: Expected exit time from a strip}

Suppose $U$ is a domain in $\CC$ which contains a point $a$. The (random) exit time $\tau(U) = \inf\{t \geq 0: B_t \notin U\}$ of a Brownian motion starting at $a$ gives some sort of measure of the size and geometry of $U$, as well as the proximity of $a$ to $U^c$. On the other hand, if $U \subsetneq \CC $ is simply connected, then the Riemann Mapping Theorem assures us of a conformal map $f$ taking $\DD = \{|z| < 1\}$ to $U$ sending $0$ to $a$. It is natural to try to use $f$ to study $\tau(U)$, and happily there is an elegant and simple connection between $f$ and $E_a[\tau(U)]$. Recall that any analytic function defined on $\DD$ admits a power series expansion $f(z) = \sum_{n=0}^{\ff} a_n z^n$ which converges for all $z$ in $\DD$. We then have the following lemma, which originally appeared in \cite[Lem. 1.1]{drum}.

\begin{lemma} \label{bigguy}
Suppose $f(z) = \sum_{n=0}^{\ff} a_n z^n$ is conformal on $\DD$. Then

\begin{equation*} \label{re2}
E_{f(0)}[\tau(f(\DD))] = \frac{1}{2}\sum_{n=1}^{\ff} |a_n|^2.
\end{equation*}

\end{lemma}

The proof in \cite{drum} uses Green's function (see Proof 4 below), but one may also obtain this result by noting that $|B_t|^2 - 2t$ is a martingale and applying the optional stopping theorem: if $0<r<1$ and $\DD_r = \{|z|<r\}$, then

\begin{equation*} \label{}
\begin{split}
E_{a}[\tau(f(\DD_r))] & = \frac{1}{2} \Big(E_{a}[|B_{\tau(f(\DD_r))}|^2]-|a|^2\Big)  \\
& = \frac{1}{2} \Big(E_{0}[|f(B_{\tau(\DD_r)})|^2]-|a|^2\Big) \\
& = \frac{1}{2} \Big( \frac{1}{2 \pi} \int_{0}^{2\pi} |f(re^{i\th})|^2 d \th -|a|^2 \Big) \\
& = \frac{1}{2} \sum_{n=1}^{\ff} |a_n|^2 r^{2n},
\end{split}
\end{equation*}

where we used Parseval's identity for analytic functions (\cite[Thm. 10.22]{rud}) as well as rotational invariance of Brownian motion (so that $B_{\tau(\DD_r)}$ is uniformly distributed on $\{|z| = r\}$) . The monotone convergence theorem allows $r$ to be taken to 1, and completes the proof.

\vski

As for Euler's sum, we let $U = \{ \frac{-\pi}{4} < \mbox{Re } z < \frac{\pi}{4} \}$. The function $\tan z$ maps $U$ conformally to $\DD$. This can be seen using the identity

\begin{equation*} \label{}
\tan z = \frac{\sin z}{\cos z} = -i \frac{e^{2iz} - 1}{e^{2iz} + 1}
\end{equation*}

and then noting that the function $x+iy \lar e^{2i(x+iy)} = e^{-2y + 2ix}$ maps $U$ conformally to $\{ \mbox{Re } z > 0\}$, and then that the M\"obius transformation $z \lar -i(\frac{z-1}{z+1})$ maps $\{ \mbox{Re } z > 0\}$ conformally to $\DD$. We therefore may apply Lemma \ref{bigguy} to $\tan^{-1}z$, which admits the Taylor series expansion

\begin{equation*} \label{}
\tan^{-1}z = z - \frac{z^3}{3} + \frac{z^5}{5} - \ldots = \sum_{n=1}^{\ff} \frac{(-1)^{n+1} z^{2n-1}}{2n-1}.
\end{equation*}

Thus,

\begin{equation*} \label{doubletrue}
E_0[\tau(U)] = \frac{1}{2} \sum_{n=1}^{\ff} \frac{1}{(2n-1)^2}.
\end{equation*}

We have arrived at a sum that is equivalent to Euler's, but we need another way to calculate $E_0[\tau(U)]$. This turns out to be easy enough: the one dimensional process $Re(B_t)^2 - t$ is a martingale, and if $T = \inf\{t \geq 0: |Re(B_t)| = \frac{\pi}{4}\}$ then we may apply the optional stopping theorem to conclude that $E_0[T] = \frac{\pi^2}{16}$. But $T$ is the same as $\tau(U)$, and we obtain

\begin{equation*} \label{}
\sum_{n=1}^{\ff} \frac{1}{(2n-1)^2} = \frac{\pi^2}{8}.
\end{equation*}

This is equivalent to Euler's sum, since

\begin{equation*} \label{}
\sum_{n=1}^\ff \frac{1}{(2n-1)^2} = \sum_{n=1}^\ff \frac{1}{n^2} - \sum_{n=1}^\ff \frac{1}{(2n)^2} = (1-\frac{1}{4})\sum_{n=1}^\ff \frac{1}{n^2} = \frac{3}{4}\sum_{n=1}^\ff \frac{1}{n^2}.
\end{equation*}

\section*{Proof 2: Exit distribution from a disk}

Analytic functions can be used to study the distribution of planar Brownian motion at certain stopping times. Suppose $\tau$ is a stopping time such that $B_\tau \in \ga$ a.s., and let $\rho^a_\tau(w) ds$ be the density of $B_\tau$ on $\ga$, when it exists, with $ds$ denoting the arclength element. Then we have the following result from \cite{medist}.

\begin{theorem} \label{massage}
Let $U$ be a domain, and suppose $f$ is a nonconstant function analytic on $U$. Let $B_t$ be a Brownian motion starting at $a$, and $\tau$ a stopping time such that the set of Brownian paths $\{B_t: 0 \leq t \leq \tau\}$ lie within $U$ a.s. Suppose that $\ga$ is a smooth curve in $U$ such that $B_\tau \in \ga$ a.s. Then, for any $a \in U$, the density for the distribution of $\hat B_{\si(\tau)} = f(B_\tau)$ on $f(\ga)$ is given by

\begin{equation} \label{cherry}
\rho^{f(a)}_{\si(\tau)} (w)ds = \sum_{z \in f^{-1}(w)\cap \ga} \frac{\rho^a_\tau (z)}{|f'(z)|} ds.
\end{equation}

\end{theorem}

The proof is essentially immediate, since any Brownian path which finishes at $f^{-1}(w) \cap \ga$ at time $\tau$ will be mapped under $f$ to a path finishing at $w$. The $|f'(z)|$ in the denominator on the right side of \rrr{cherry} is the scaling factor required for the change in the arclength element mapped under the analytic function $f$. If $f$ is conformal then it is often easier to use $g = f^{-1}$, and \rrr{cherry} becomes

\begin{equation*} \label{cherry2}
\rho^{a}_{\si(\tau)} (w)ds = \rho^{g(a)}_{\tau} (g(w))|g'(w)| ds.
\end{equation*}

Armed with this, we can find the density of the distribution of Brownian motion at a large number of stopping times. There is only one which is obvious, namely $\rho^0_{\tau(\DD)}(e^{i \th})ds = \frac{1}{2\pi} ds$, and this is evident by the rotational invariance of Brownian motion. This may now be used along with the formula $\psi(z) = \frac{z-a}{1-\bar a z}$ for a conformal self-map of $\DD$ sending $a$ to 0 in order to derive

\begin{equation} \label{mae}
\rho^a_{\tau(\DD)}(e^{i \th})ds = \frac{1}{2\pi} \frac{1-|a|^2}{|1-\bar a e^{i\th}|^2} ds.
\end{equation}

Let us now find the density for the exit time of $\HH = \{Im(z)>0\}$. The conformal map taking $\HH$ to $\DD$ is given by $f(z) = \frac{z-i}{z+i}$, with $f'(z) = \frac{2i}{(z+i)^2}$. We obtain

\begin{equation*} \label{}
\rho^i_{T_\HH}(x)ds = \frac{1}{\pi} \frac{1}{1+x^2}ds.
\end{equation*}

To find the distribution from an arbitrary point on the positive imaginary axis $vi$ use the map $f(z) = vz$ to obtain

\begin{equation} \label{swimchick}
\rho^{vi}_{T_\HH}(x)ds = \frac{1}{\pi} \frac{v}{v^2+x^2}ds.
\end{equation}

Now we will use this in order to calculate the density for the exit time of the punctured disk $\DD^\times = \{0 < |z| < 1\}$ when the Brownian motion is started at a point $a \in (0,1)$. If we take $v = - \ln a$, then $f(z) = e^{iz}$ maps $\HH$ onto $\DD^\times$, wrapping $\dd \HH$ around $\dd \DD$ and mapping $vi$ to $a$. The preimages of a point $e^{i \th}$ under $f$ are all points of the form $\th+2\pi k$ for integer $k$, and $|f'| = 1$ on $\dd \HH$. Theorem \ref{massage} therefore gives

\begin{equation} \label{page}
\rho^a_{\tau(\DD^\times)}(e^{i \th})ds= \sum_{k=-\ff}^{\ff} \frac{-\ln a}{\pi ((\ln a)^2 + (\th+2\pi k)^2)}ds.
\end{equation}

On the other hand, we may calculate $\rho^a_{\tau(\DD^\times)}(e^{i \th})ds$ in a different way, by noting that for $a \neq 0$ it coincides with the exit distribution of the disk $\DD$, as $P_a(B_{\tau(\DD^\times)} = 0) = 0$ (see \cite[Sec. 5.1]{durBM}). The quantities \rrr{mae} and \rrr{page} must therefore be equal, the difference in form due to the fact that each term in the sum in \rrr{page} corresponds to a different homotopy class of paths in the punctured disk terminating at $e^{i\th}$, while \rrr{mae} does not differentiate between the homotopy classes. Equating them, assuming $a \in (0,1)$ for simplicity, gives the identity

\begin{equation*} \label{effyou}
\sum_{k=-\ff}^{\ff} \frac{-\ln a}{\pi ((\ln a)^2 + (\th+2\pi k)^2)} = \frac{1}{2\pi} \frac{1-a^2}{|1-a e^{i\th}|^2} = \frac{1-a^2}{2\pi(1+a^2-2a\cos \th)}.
\end{equation*}

Divide both sides by $-\ln a$ and simplify. This gives

\begin{equation*} \label{cthg}
\sum_{k=-\ff}^{\ff} \frac{1}{((\ln a)^2 + (\th + 2\pi k)^2) } = \frac{1-a^2}{2(-\ln a)(1+a^2-2a\cos \th)}.
\end{equation*}

Assuming $\th \neq 0$, we can now let $a \nearrow 1$ using $\lim_{a \lar 1} \frac{1-a^2}{\ln a} = -2$, and obtain

\begin{equation*} \label{cthg2}
\sum_{k=-\ff}^{\ff} \frac{1}{(\th + 2\pi k)^2 } = \frac{1}{2(1-\cos \th)}.
\end{equation*}

Subtract $\frac{1}{\th^2}$ from both sides (the term $k=0$), take the limit as $\th \lar 0$ using $\lim_{\th \lar 0} \frac{1}{2(1-\cos \th)} - \frac{1}{\th^2} = \frac{1}{12}$, and simplify. We obtain

\begin{equation*} \label{jill}
\sum_{k=1}^{\ff} \frac{1}{k^2} = \frac{\pi^2}{6}.
\end{equation*}

\section*{Proof 3: Exit distribution from a strip}

We now calculate the exit distribution on an infinite strip. We let $W = \{-1 < Re(z) < 1\}$ and take the starting point of the Brownian motion $a$ to lie on the real interval $(-1,1)$. We will calculate density of the exit distribution at the boundary point 1 by applying our theorem to the function $\tan (\frac{\pi}{4} z)$, which maps $W$ conformally to $\DD$, to get

\begin{equation} \label{property}
\begin{split}
\rho^a_{\tau(W)}(1)ds & = \frac{\pi}{4}\sec^2(\frac{\pi}{4})\rho^{\tan (\frac{\pi}{4} a)}_{\tau(\DD)} (\tan (\frac{\pi}{4}))ds \\
& = \frac{\pi}{2} \rho^{\tan (\frac{\pi}{4} a)}_{\tau(\DD)} (1) ds \\
& = \frac{1-\tan^2(\frac{\pi}{4} a)}{4(1-\tan (\frac{\pi}{4} a))^2}ds \\
& = \frac{1+\tan(\frac{\pi}{4} a)}{4(1-\tan (\frac{\pi}{4} a))}ds,
\end{split}
\end{equation}

where the exit distribution of the disk given in Proof 2 was used. There is another method for calculating $\rho^a_{T_W}(\pm 1+yi)ds$ which uses a form of the reflection principle. For real number $b$ let us define $\tau(b) = \inf\{ t \geq 0 : Re(B_t) = b\}$, and then recursively define $\tau(b_1, b_2, \ldots, b_{n+1}) = \inf\{t \geq \tau(b_1, \ldots, b_{n}): Re(B_t) = b_{n+1}\}$; that is, $\tau(b_1, b_2, \ldots, b_{n})$ is the first time at which $Re(B_t)$ has visited the sequence $b_1, b_2, \ldots , b_n$ in order. We then claim

\begin{equation} \label{denim}
\rho^a_{T_W}(1)ds = \rho^a_{\tau(1)}(1)ds - \rho^a_{\tau(-1,1)}(1)ds + \rho^a_{\tau(1,-1,1)}(1)ds - \rho^a_{\tau(-1,1,-1,1)}(1)ds + \ldots.
\end{equation}

Intuitively, this is very simple: in order to calculate $\rho^a_{T_W}(1)ds$ we want to count paths which leave $\{Re(z)<1\}$ at the point $1$, however we want to remove the contribution from paths which strike $\{Re(z)=-1\}$ first, so we consider $\rho^a_{\tau(1)}(1)ds - \rho^a_{\tau(-1,1)}(1)ds$; however, by subtracting $\rho^a_{\tau(-1,1)}(1)ds$ we have subtracted too much, as we have incorrectly subtracted the contribution from paths which hit $\{Re(z)=1\}$ before $\{Re(z)=-1\}$, so we must add $\rho^a_{\tau(1,-1,1)}(1)ds$ to compensate; however, by an analogous argument we have overcompensated, and must therefore subtract $\rho^a_{\tau(-1,1,-1,1)}(1)ds$, and so forth. The reader with a penchant for rigor can find a complete proof in \cite{medist}. Now, by the reflection principle (see \cite{fima} or \cite{mortper}), the process reflected over $\{Re(z) = 1\}$ after time $\tau(1)$, i.e.

\begin{equation*}
\hat B_{t} = \left \{ \begin{array}{ll}
B_t & \qquad  \mbox{if } t \leq \tau(1)  \\
1+(1-Re(B_t)) + i Im(B_t) & \qquad \mbox{if } t > \tau(1)\;,
\end{array} \right.
\end{equation*}

is also a Brownian motion, and it is therefore evident that the distribution of $\hat B_{\tau(1,3)} = \hat B_{\tau(3)}$ is identical to that of $\hat B_{\tau(1,-1)}$. In this manner, by reflection, the quantities on the right side of \rrr{denim} can all be expressed in terms of the exit distribution of a half plane, and we obtain

\begin{equation} \label{tats}
\rho^a_{T_W}(1)ds = \rho^a_{\tau(1)}(1)ds - \rho^a_{\tau(-3)}(-3)ds + \rho^a_{\tau(5)}(5)ds - \rho^a_{\tau(-7)}(-7)ds + \ldots
\end{equation}

We need to calculate $\rho^a_{\tau(r)}(r)ds$ for any real $r$, and this was done already in Proof 2: it is immediate from \rrr{swimchick} by rotation and translation that $\rho^a_{\tau(r)}(r)ds = \frac{ds}{\pi |r-a|}$. Plugging this into \rrr{tats} gives

\begin{equation*} \label{inside}
\rho^a_{T_W}(1)ds = \frac{ds}{\pi} \Big( \frac{1}{1-a} - \frac{1}{3+a} + \frac{1}{5-a} - \frac{1}{7+a} + \ldots \Big)
\end{equation*}

Equating this with \rrr{property} yields

\begin{equation} \label{mei}
\sum_{j=1}^{\ff} \frac{(-1)^{j+1}}{(2j-1) + (-1)^j a } = \frac{1}{1-a} - \frac{1}{3+a} + \frac{1}{5-a} - \frac{1}{7+a} + \ldots=  \frac{\pi}{4} \Big(\frac{1+\tan (\frac{\pi}{4} a)}{1-\tan (\frac{\pi}{4} a)}\Big).
\end{equation}

We now notice that taking $a=0$ immediately gives us Leibniz's representation for $\pi$:

\begin{equation*} \label{propasi}
\frac{\pi}{4} = 1-\frac{1}{3} + \frac{1}{5} - \frac{1}{7} + \ldots
\end{equation*}

However, this sum, interesting though it may be, is not the one we are after. In order to obtain Euler's sum, we simply differentiate both sides of \rrr{mei} with respect to $a$ to obtain

\begin{equation*} \label{mei3}
\sum_{j=1}^{\ff} \frac{1}{((2j-1) + (-1)^j a)^2} = \frac{1}{(1-a)^2} + \frac{1}{(3+a)^2} + \frac{1}{(5-a)^2} + \frac{1}{(7+a)^2} + \ldots=  \frac{\pi^2}{8} \Big(\frac{\sec^2(\frac{\pi}{4} a))}{(1-\tan (\frac{\pi}{4} a))^2}\Big).
\end{equation*}

Setting $a=0$ at this point now gives the sum of the odd squares, which was shown at the end of Proof 1 to be equivalent to Euler's sum. In fact, it should be clear that successive differentiation will lead to the values of all sums of the form $\sum_{j=1}^{\ff} \frac{1}{(2j-1)^{2m}}$ or $\sum_{j=1}^{\ff} \frac{(-1)^j}{(2j-1)^{2m+1}}$ for integer $m$.

\section*{Proof 4: Green's function of a disk}

This proof is the spiritual cousin of Proof 2. As before, let us run a Brownian motion $B_t$ until a stopping time $\tau$, and let $B_t = \Delta$ for $t \geq \tau$, where $\Delta$ is a so-called "cemetery point" outside of $\CC$; in other words, $B_t$ is no longer in the plane for $t \geq \tau$. Let $\rho_t^\tau(z,w)$ be the probability density function at point $w$ and time $t$ of the distribution of this killed Brownian motion. We may define the {\it occupation measure} of $B_t$ up to time $\tau$ as $\mu(A) = \int_{0}^{\tau} 1_A(B_s) ds$; that is, $\mu$ is a random measure representing the amount of time $B_t$ spends in $A$ up until time $\tau$. Calculating formally,

\begin{equation*} \label{exploctime}
\begin{split}
E_a \int_{0}^{\tau} 1_A(B_s) ds & =  \int_{0}^{\tau} E_a[1_A(B_s)] ds \\
& = \int_{0}^{\ff} \int_{\CC} 1_A (z) \rho_s^\tau(a,z) dA(z) ds \\
& = \int_{\CC} \Big( \int_{0}^{\ff} \rho_s^\tau(a,z) ds \Big) 1_A(z) dA(z),
\end{split}
\end{equation*}

where $dA(z)$ signifies an area integral in the plane. It follows that the function

\begin{equation} \label{probdef}
G_\tau(a,z) := \int_{0}^{\ff} \rho_s^\tau(a,z) ds
\end{equation}

is the density of the measure $E_a[\mu(A)]$. We call \rrr{probdef} the {\it Green's function} associated with $\tau$; the name is due to the fact that when the stopping time in question is the exit time of a domain, \rrr{probdef} agrees up to a multiplicative constant with the Green's function of that domain. Green's function is well known among analysts to be conformally invariant, however associating it to a stopping time rather than a domain allows us to deduce more than this from L\'evy's Theorem. The following is a simplification of \cite[Prop. 2]{megreen}.

\begin{theorem} \label{chimass}
Let $\Om$ be a domain, and suppose $f$ is a function analytic on $\Om$ with $f'(z) \neq 0$ for all $z \in \Om$. Let $B_t$ be a Brownian motion starting at $a$, and $\tau$ a stopping time such that the set of Brownian paths $\{B_t: 0 \leq t \leq \tau\}$ lie within $U$ a.s. Then

\begin{equation*} \label{cherry2}
G_{\si(\tau)} (f(a),w) = \sum_{w' \in f^{-1}(\{w\})} G_\tau (a,w').
\end{equation*}

\end{theorem}

See \cite{megreen} for the straightforward proof of this result. It is clear that this result can allow us to calculate a large number of Green's functions, but, as was the case when we considered the distribution at stopping times in Proof 2, we need a simple one to begin with. In that case we used the exit distribution from a disk, but in this case it is easier to use a half-plane. Recall that $\HH = \{Im(z)>0\}$. The density $\rho_t^{\tau(\HH)} (a,z)$ must capture the probability of Brownian paths near $z$ at time $t$, but only those which have not previously intersected $\RR$. By the reflection principle we used in Proof 3, the processes $B_t$ and

\begin{equation*}
\hat B_t = \left \{ \begin{array}{ll}
B_t & \qquad  \mbox{if } t \leq {\tau(\HH)}  \\
\overline B_t & \qquad \mbox{if } t > {\tau(\HH)} \;,
\end{array} \right.
\end{equation*}

have the same law, where $\overline B_t$ is simply the complex conjugation map $x + yi \lar x - yi$ applied to the complex-valued process $B_t$. But $B_t$ is near $z$ and $t > \tau(\HH)$ precisely when $\hat B_t$ is near $\bar z$ with $t > \tau(\HH)$, and this occurs precisely when $\hat B_t$ is near $\bar z$, since the Brownian motion cannot travel from $a$ to $\bar z$ without first crossing $\RR$. We conclude that $\rho_t^{\tau(\HH)} (a,z) = \rho_t (a,z) - \rho_t (a,\bar z) = \frac{1}{2\pi t} (e^{-|a-z|^2/(2t)} - e^{-|a-\bar z|^2/(2t)})$. We therefore have

\begin{equation} \label{hp}
\begin{split}
G_{\tau(\HH)}(a,z) & = \int_{0}^{\ff} \frac{1}{2\pi t} (e^{\frac{-|a-z|^2}{2t}} - e^{\frac{-|a-\bar z|^2}{2t}}) dt \\
& = \frac{1}{2\pi} \int_{0}^{\ff} \frac{1}{t} \int_{|a-z|}^{|a-\bar z|} \frac{s}{t}e^{\frac{-s^2}{2t}} ds dt \\
& = \frac{1}{2\pi} \int_{|a-z|}^{|a-\bar z|} \int_{0}^{\ff} \frac{s}{t^2} e^{\frac{-s^2}{2t}} dt ds \\
& = \frac{1}{\pi} \int_{|a-z|}^{|a-\bar z|} \frac{1}{s} \int_{0}^{\ff} e^{-u} du ds \\
& = \frac{1}{\pi} \ln \Big( \frac{|a-\bar z|}{|a-z|}\Big).
\end{split}
\end{equation}

$\phi(z) = -i \Big(\frac{z-1}{z+1}\Big)$ maps $\DD$ conformally onto $\HH$, sending $0$ to $i$. Thus,

\begin{equation*} \label{}
\begin{split}
G_{\tau(\DD)}(0,z) & =G_{\tau(\HH)}(i,\phi(z)) = \frac{1}{\pi} \ln \Big( \frac{|i-i(\frac{\bar z -1}{\bar z + 1})|}{|i+i(\frac{z -1}{z + 1})|}\Big) \\
& = \frac{1}{\pi} \ln \Big| \frac{2z+2}{2|z|^2 + 2z}\Big| = \frac{1}{\pi} \ln \frac{1}{|z|}.
\end{split}
\end{equation*}

To calculate $G_{\tau(\DD)}(a,w)$ for arbitrary $a \in \DD$, we use the disk automorphism $\phi(z) = \frac{z-a}{1-\bar a z}$ we used in Proof 2. It follows that

\begin{equation*} \label{lamb}
G_{\tau(\DD)}(a,z) = G_{\tau(\DD)}(\phi(a),\phi(z)) = \frac{1}{\pi} \ln \frac{|1-\bar a z|}{|z-a|}.
\end{equation*}

Now, for $a \neq 0$ and $\DD^\times = \{0 < |z| < 1\}$ we have $G_{\tau(\DD)}(a,z) = G_{\tau(\DD^\times)}(a,z)$, for the same reason as was used in Proof 2: $P_a(B_{\tau(\DD^\times)} = 0) = 0$. Furthermore, we can again use the covering map $f(z) = e^{iz}$ in order to calculate $G_{\tau(\DD^\times)}(a,z)$. Let us assume for simplicity that $a,z \in (0,1)$, so we may write $a=e^{-\aaa}, z = e^{-\ga}$, for $\aaa, \ga > 0$. We can take $\aaa i$ as the preimage of $a$ under $f$, and $\ga i + 2\pi n$ for integer $n$ will therefore be all preimages of $z$ under $f$. Using \rrr{hp} and Theorem \ref{chimass} we have

\begin{equation*} \label{}
\begin{split}
G_{\tau(\DD^\times)}(e^{-\aaa},e^{-\ga}) & = \frac{1}{\pi} \sum_{n=-\ff}^{\ff} \ln \frac{|(\aaa+\ga)i - (2\pi n)|}{|(\aaa-\ga)i - (2\pi n)|} \\
& = \frac{1}{2\pi} \sum_{n=-\ff}^{\ff} \ln \frac{(2\pi n)^2 + (\aaa+\ga)^2}{(2\pi n)^2 + (\aaa-\ga)^2}.
\end{split}
\end{equation*}

This quantity must therefore be equal to $\frac{1}{\pi} \ln \frac{|1-e^{-(\aaa+\ga)}|}{|e^{-\ga}-e^{-\aaa}|}$. Multiplying each of these by $2\pi$ and exponentiating gives the identity

\begin{equation*} \label{mirror}
\prod_{n=-\ff}^\ff \frac{(2\pi n)^2+(\aaa+\ga)^2}{(2\pi n)^2+(\aaa-\ga)^2} = \Big(\frac{1-e^{-(\aaa+\ga)}}{e^{-\aaa} - e^{-\ga}}\Big)^2.
\end{equation*}

Now multiply top and bottom of the expression on the right by $(e^{\frac{(\aaa+\ga)}{2}})^2$, and divide top and bottom of all terms on the left other than $n=0$ by $(2 \pi n)^2$. We obtain

\begin{equation*} \label{}
\frac{(\aaa+\ga)^2}{(\aaa-\ga)^2} \Big(\prod_{n=1}^\ff \frac{(1+(\frac{\aaa+\ga}{2\pi n})^2)}{(1+(\frac{\aaa-\ga}{2\pi n})^2)}\Big)^2 = \frac{\sinh (\frac{\aaa+\ga}{2})^2}{\sinh (\frac{\aaa-\ga}{2})^2}.
\end{equation*}

Multiplying both sides by $(\aaa-\ga)^2$ and taking the limit as $\ga \lar \aaa$ leads to the infinite product representation for $\sinh$:

\begin{equation} \label{chinchin}
\sinh \aaa = \aaa \prod_{n=1}^\ff (1+(\frac{\aaa}{\pi n})^2).
\end{equation}

Note that from this, using $\sin \aaa = i \sinh(-i\aaa)$, we can deduce \rrr{sineprod}, the infinite product representation for sine, and thus Euler's sum. However, Euler's sum follows just as easily from \rrr{chinchin} directly, for we may write $\frac{\sinh \aaa}{\aaa} = 1+ \frac{\aaa^2}{3!} + \frac{\aaa^4}{5!} + \ldots$, while the coefficient of $\aaa^2$ in the expansion of the infinite product in \rrr{chinchin} is $\sum_{n=1}^{\ff}\frac{1}{(\pi n)^2}$.


\section*{Other proofs and concluding remarks}

There are several other proofs of Euler's sum using Brownian motion available. For example, calculating the exit distribution of $\{-1<Re(z)<1, Im(z) > 0\}$ on the boundary component $\{-1<Re(z)<1, Im(z) = 0\}$ using reflection as in Proof 3 leads to an identity taking Euler's sum as a special case, as does calculating the exit distribution of $\CC \backslash [-1,1]$ using a covering map as in Proof 2. Details can be found in \cite{medist}. Doubtless there are others using exit distributions.

\vski

It is possible to calculate the Green's function for a strip using reflection in a similar manner to what was done in Proof 3. As in Proof 4, this leads to an infinite product identity, but it is not clear whether Euler's identity lives inside this product. One difficulty is that the product is not absolutely convergent, and care must be taken in passing to the limit; this was not adequately explained in \cite{megreen}, but is discussed in detail in \cite{melgreen} and \cite{melgreen2}, where analytic methods rather than probabilistic ones are used to obtain the same result. Incidentally, it would be interesting to see whether the analytic methods in \cite{melgreen} and \cite{melgreen2} can be adapted to mimic Proof 4 above; the difficulty seems to lie in coping with the center point 0 in the disk, while using Brownian motion allows us to disregard this point.

\vski

The author has been asked on more than one occasion why Euler's sum in particular keeps popping up when working with planar Brownian motion. The short and correct answer is "I don't know": when you play around with Brownian motion in the complex plane looking for identities, you don't really get to choose what comes out. However, in a sense, it is natural: the Hilbert-space techniques inherent in the theory ensure that quantities have a tendency to become squared when working with Brownian motion, so if you can bring the integers into play you have a chance of happening upon Euler's sum. Proof 1 is, after all, a Fourier series proof in disguise, since that's what complex power series are; Proofs 2 and 4 bring the integers into the mix due to homotopy considerations, the winding number of Brownian motion about the origin; and the integers appear in Proof 3 in regards to the number of times the Brownian motion has been reflected over the lines $\{Re(z) = \pm 1\}$.


\section*{Acknowledgements}

I am grateful for support from Australian Research Council Grants DP0988483 and DE140101201. I would also like to express my thanks to Maher Boudabra, Paul Bourgade, and Jim Pitman for helpful comments.

\bibliographystyle{alpha}
\bibliography{CAbib}

\providecommand{\bysame}{\leavevmode\hbox to3em{\hrulefill}\thinspace}
\providecommand{\MR}{\relax\ifhmode\unskip\space\fi MR }
\providecommand{\MRhref}[2]{%
  \href{http://www.ams.org/mathscinet-getitem?mr=#1}{#2}
}
\providecommand{\href}[2]{#2}
\begin{thebibliography}{10}

\bibitem{drum}
R.~Ba{\~n}uelos and T.~Carroll, \emph{Brownian motion and the fundamental
  frequency of a drum}, Duke Mathematical Journal \textbf{75} (1994), no.~3,
  575--602.

\bibitem{bassy}
R.~Bass, \emph{Probabilistic techniques in analysis}, Springer-Verlag, 1995.

\bibitem{biane}
P.~Biane, J.~Pitman, and M.~Yor, \emph{Probability laws related to the {J}acobi
  theta and {R}iemann zeta functions, and {B}rownian excursions}, Bulletin of
  the American Mathematical Society \textbf{38} (2001), no.~4, 435--465.

\bibitem{bourg}
P.~Bourgade, T.~Fujita, and M.~Yor, \emph{Euler's formulae for {$\zeta(2n)$}
  and products of {C}auchy variables}, Electronic Communications in Probability
  \textbf{12} (2007), 73--80. \MR{2300217}

\bibitem{chappy}
R.~Chapman, \emph{Evaluating $\zeta (2)$}, available at
  http://secamlocal.ex.ac.uk/people/staff/rjchapma/etc/zeta2.pdf.

\bibitem{durBM}
R.~Durrett, \emph{Brownian motion and martingales in analysis}, Wadsworth
  Advanced Books \& Software, 1984.

\bibitem{yano}
T.~Fujita and Y.~Yano, \emph{Special values of the {H}urwitz zeta function via
  generalized {C}auchy variables}, Kyoto Journal of Mathematics \textbf{52}
  (2012), no.~3, 465--477.

\bibitem{fima}
F.~Klebaner, \emph{Introduction to stochastic calculus with applications},
  Imperial College Press, 2005.

\bibitem{medist}
G.~Markowsky, \emph{On the distribution of planar {B}rownian motion at stopping
  times}, Annales Academi\ae { }Scientiarum Fennic\ae { }Mathematica, to
  appear, arxiv/1606.03186.

\bibitem{meecp}
G.~Markowsky, \emph{On the expected exit time of planar {B}rownian motion from
  simply connected domains}, Electronic Communications in Probability
  \textbf{16} (2011), 652--663.

\bibitem{megreen}
G.~Markowsky, \emph{On the planar {B}rownian {G}reen's function for stopping
  times}, Journal of Mathematical Analysis and Applications \textbf{455}
  (2017), no.~2, 1221 -- 1233.

\bibitem{melgreen}
Y.~Melnikov, \emph{A new approach to the representation of trigonometric and
  hyperbolic functions by infinite products}, Journal of Mathematical Analysis
  and Applications \textbf{344} (2008), no.~1, 521--534.

\bibitem{melgreen2}
{Y}. Melnikov, \emph{Green's functions and infinite products: bridging the
  divide}, Springer Science \& Business Media, 2011.

\bibitem{mortper}
P.~M{\"o}rters and Y.~Peres, \emph{Brownian motion}, vol.~30, Cambridge
  University Press, 2010.

\bibitem{pace}
L.~Pace, \emph{Probabilistically proving that $\zeta(2) = \pi^2/6$}, American
  Mathematical Monthly \textbf{118} (2011), no.~7, 641--643.

\bibitem{pityor}
J.~Pitman and M.~Yor, \emph{Infinitely divisible laws associated with
  hyperbolic functions}, Canadian Journal of Mathematics \textbf{55} (2003),
  no.~2, 292--330.

\bibitem{revyor}
D.~Revuz and M.~Yor, \emph{{Continuous martingales and Brownian motion}},
  Springer Verlag, 1999.

\bibitem{rud}
W.~Rudin, \emph{Real and complex analysis}, Tata McGraw-Hill, 2006.

\bibitem{schaum}
M.~Spiegel, \emph{Complex variables: {W}ith an introduction to conformal
  mapping and its applications}, Schaum's Outline Series, 2009.

\bibitem{stein}
E.~Stein and R.~Shakarchi, \emph{Fourier analysis: an introduction}, vol.~1,
  Princeton University Press, 2011.

\end{thebibliography}

\end{document}